\documentclass[11pt]{article}
\usepackage{makeidx}
\usepackage{epsfig}
\usepackage{float}
\usepackage{amscd}
\usepackage{amsmath}
\usepackage{amssymb}
\usepackage{amsthm}
\usepackage{mathtools}
\usepackage{color}
\usepackage{fancyheadings}
\usepackage{fancyvrb}
\usepackage{diagbox}
\usepackage{url}

\newcommand{\omitit}[1]{}

\newcommand{\norm}[1]{\Vert{#1}\Vert}

\newcommand{\rr}{{\mathbf r}}
\newcommand{\ww}{{\mathbf w}}
\newcommand{\ba}{{\mathbf a}}
\newcommand{\xx}{{\mathbf x}}
\newcommand{\yy}{{\mathbf y}}
\newcommand{\cK}{{\mathcal K}}
\DeclareMathOperator{\spa}{span}
\DeclareMathOperator{\proj}{proj}
\DeclareMathOperator{\diag}{diag}

\numberwithin{equation}{section}
\newtheorem{theorem}{Theorem}[section]
\newtheorem{alg}[theorem]{Algorithm}

\date{}

\makeindex
\usepackage{makeidx}

\begin{document}
\title{Extrapolating the Arnoldi Algorithm \\
To Improve Eigenvector Convergence}
\author{Sara Pollock\\Department of Mathematics, University of Florida\\
\\
L. Ridgway Scott\\ University of Chicago}
 \def\thepage{}\maketitle\pagenumbering{arabic}
\centerline{\today}


\abstract{We consider extrapolation of the Arnoldi algorithm to accelerate computation of the dominant eigenvalue/eigenvector pair.
The basic algorithm uses sequences of Krylov vectors to form a small eigenproblem 
which is solved exactly.
The two dominant eigenvectors output from consecutive Arnoldi steps are then 
recombined to form an extrapolated iterate, and this accelerated iterate is used to 
restart the next Arnoldi process.
We present numerical results testing the algorithm on a variety of cases and find on most
examples it substantially improves the performance of restarted Arnoldi.
The extrapolation is a simple post-processing step which has minimal computational cost.}


\maketitle

\section{Introduction}\label{sec:intro}
There are many applications in which the smallest eigenvalues of large systems 
must be computed, e.g., in stability analysis of numerical schemes \cite{lrsBIBiu}
and in stability analysis of partial differential equations \cite{lrsBIBis}.
The (inverse) power method is often preferred due to 
its ease of implementation and the limited amount of storage required.
In many cases \cite{lrsBIBis}, all that is required to implement the inverse power method 
is to solve the associated system of equations repeatedly.
Thus it is easy to modify a code for a solver to become a code for the inverse power method.

Extrapolation has been shown \cite{ref:Pollackextraeig} to provide an effective way to improve
the power method.
Here we consider a different approach in which the power method is first generalized.
We then examine extrapolation of the method.

One generalization of the power method is to use a small number $k>1$ of approximation
vectors and to project the eigenproblem onto the corresponding $k$-dimensional space.
One then solves the projected $k$-dimensional eigenproblem and extracts 
the eigenvector corresponding to the extreme eigenvalue as the next iterate.
This can lead to faster convergence with a controlled increase in storage.
A natural set of vectors to use is a Krylov basis,
and we dub this approach the $k$-step Krylov method.
We show that the popular LOBPCG method is of this form.
There could of course be other ways of generating appropriate vectors at each step,
e.g., a random process.

We show that the Arnoldi algorithm provides a very stable implementation of 
the $k$-step Krylov method.
We further demonstrate how a simple extrapolation technique,
which takes a combination of the two latest Arnoldi outputs as the next approximation, 
can be used to further enhance the rate of convergence.
Since the basic $k$-step method is Arnoldi, it is remarkable that this algorithm
can be improved by extrapolation.

\section{$k$-step Krylov methods}
\label{sec:kaystepm}

We can define general $k$-step Krylov methods as follows.
We start with a vector $\yy_1$ and define $\yy_{j+1}=A\yy_j$, for $j=1,\dots,k-1$.
We seek to find coefficients $a_1,\dots,a_k$ such that
\begin{equation}\label{eqn:analogpameig}
A\sum_j a_j \yy_j =\lambda \sum_j a_j \yy_j .
\end{equation}
Taking dot products, we see that this is equivalent to
$$
K\ba=\lambda M\ba,
$$
where
$$
K_{ij}=\yy_i^t A \yy_j ,\qquad M_{ij}=\yy_i^t \yy_j .
$$
Thus we can determine $\ba$ by solving the eigenproblem
\begin{equation}\label{eqn:naivempameig}
M^{-1}K\ba=\lambda \ba.
\end{equation}
Due to the close connection with Ritz methods, $\ba$ is called the Ritz
vector and $\lambda$ is the Ritz value.
A great deal is known about how these approximate eigenvectors and eigenvalues
for $A$ as $k$ increases \cite{saad1980rates}.

We define the output of the $k$-step method as $\lambda_1,\dots\lambda_m$ 
($m<k$) and
$$
 \yy=\sum_{i=1}^k a_i \yy_i,
$$
where $\ba$ is the eigenvector corresponding to the extreme eigenvalue $\lambda_1$
for \eqref{eqn:naivempameig}, and $\lambda_2,\dots\lambda_m$ are the remaining
eigenvalues in descending order.

The eigenproblem \eqref{eqn:naivempameig} can be solved analytically for $k\leq 4$.
In \cite{lrsBIBjf}, the case $k=2$ is described in detail.

\subsection{Orthogonalization of Krylov vectors}
\label{sec:orthokyl}

Unfortunately, the naive approach \eqref{eqn:naivempameig} to the $k$-step method
fails for larger $k$, due to ill conditioning, as described in \cite{lrsBIBjf}.
Thus we consider orthogonalization of the Krylov vectors.

Now we modify the steps leading to \eqref{eqn:naivempameig}.
We start with a vector $\hat \yy_1$ and define Krylov vectors
$\hat \yy_{j+1}=A\hat \yy_j$, for $j=1,\dots,k$.
Then we orthogonalize to get $\yy_1,\dots,\yy_k$ by the modified Gram--Schmidt algorithm
\cite{bjorck1994numerics,paige2006modified}.
We provide the details in \cite{lrsBIBjf}. 
In this setting, the matrix $M$ is the identity.

The use of orthogonal Krylov vectors allows extension to more steps $k$, 
but for slightly larger $k$ the algorithm still
fails due to the increasing condition number of $K$, as indicated in \cite{lrsBIBjf}. 

\subsection{Arnoldi algorithm}

The Arnoldi algorithm makes a small change in the order of orthogonalization and 
multiplication by the matrix $A$.
Instead of first creating the Krylov vectors all at once, we multiply by $A$
only after orthogonalization.
Thus $\yy_1=\norm{\hat \yy_1}^{-1}\hat \yy_1$.
Then for $n= 1,2,\dots,k-1$, define
\begin{equation}\label{eqn:larnwhydefoag}
\begin{split}
\widetilde \yy_n=A \yy_{n} - \sum_{j=1}^{n} h_{j,n} \yy_j\,, \qquad
&h_{j,n}= \yy_j^{\;t} A\yy_{n}\,,\quad j=1,\dots n, \\
h_{n+1,n}=\norm{\widetilde \yy_n},\qquad
&\yy_{n+1}=h_{n+1,n}^{-1} \widetilde \yy_n.
\end{split}
\end{equation}
For $n=k$, we compute
$h_{j,k}= \yy_j^{\;t}A\yy_{k}$ for $j=1,\dots k$,
but we do not perform the orthogonalization steps in the first line 
of \eqref{eqn:larnwhydefoag} for $n=k$.

\begin{figure}
\centerline{(a)\includegraphics[trim = 0pt 0pt 30pt 0pt,clip = true,width=0.42\textwidth]
{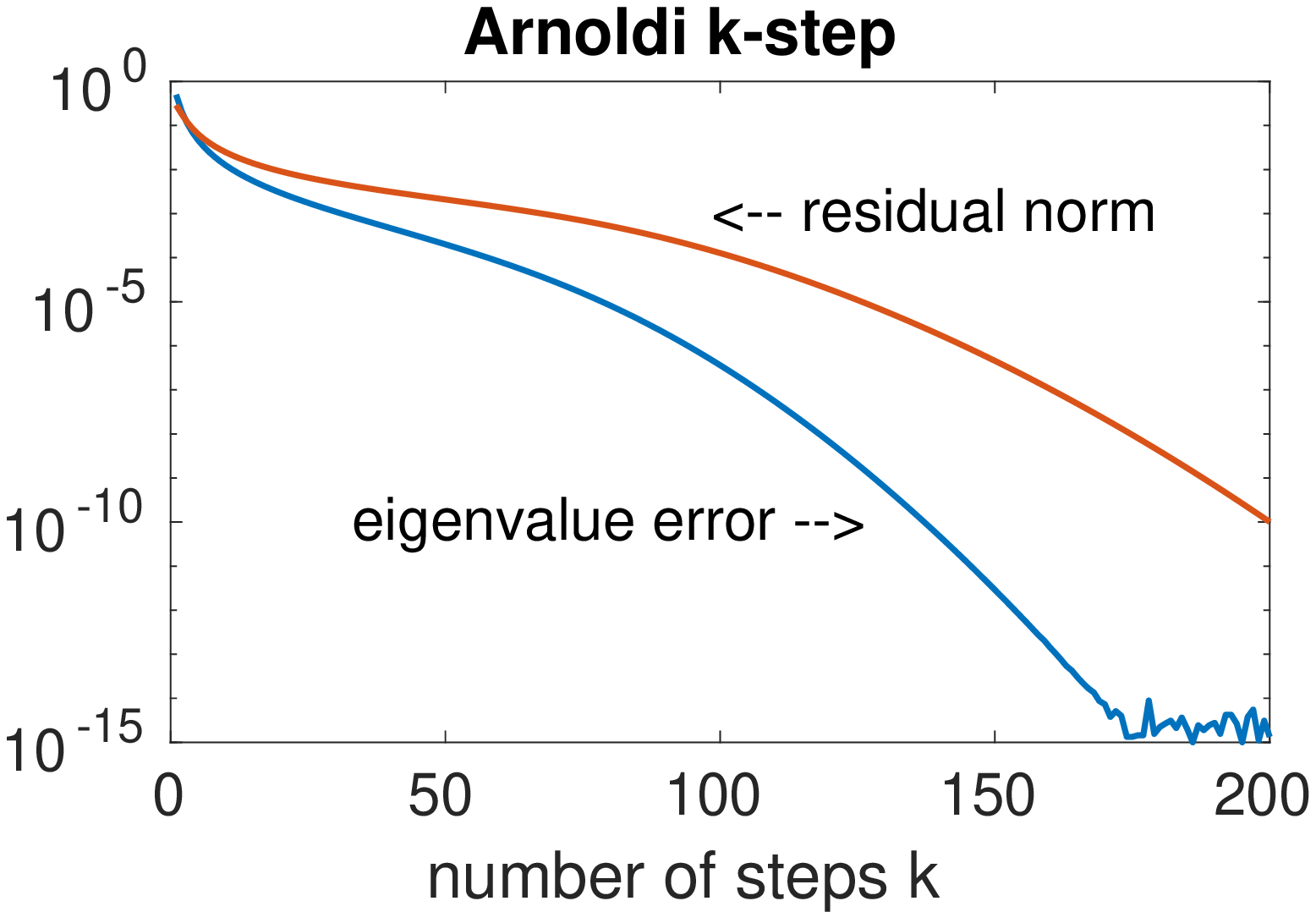}
 \qquad     
(b)\includegraphics[trim = 0pt 0pt 30pt 0pt,clip = true,width=0.42\textwidth]
{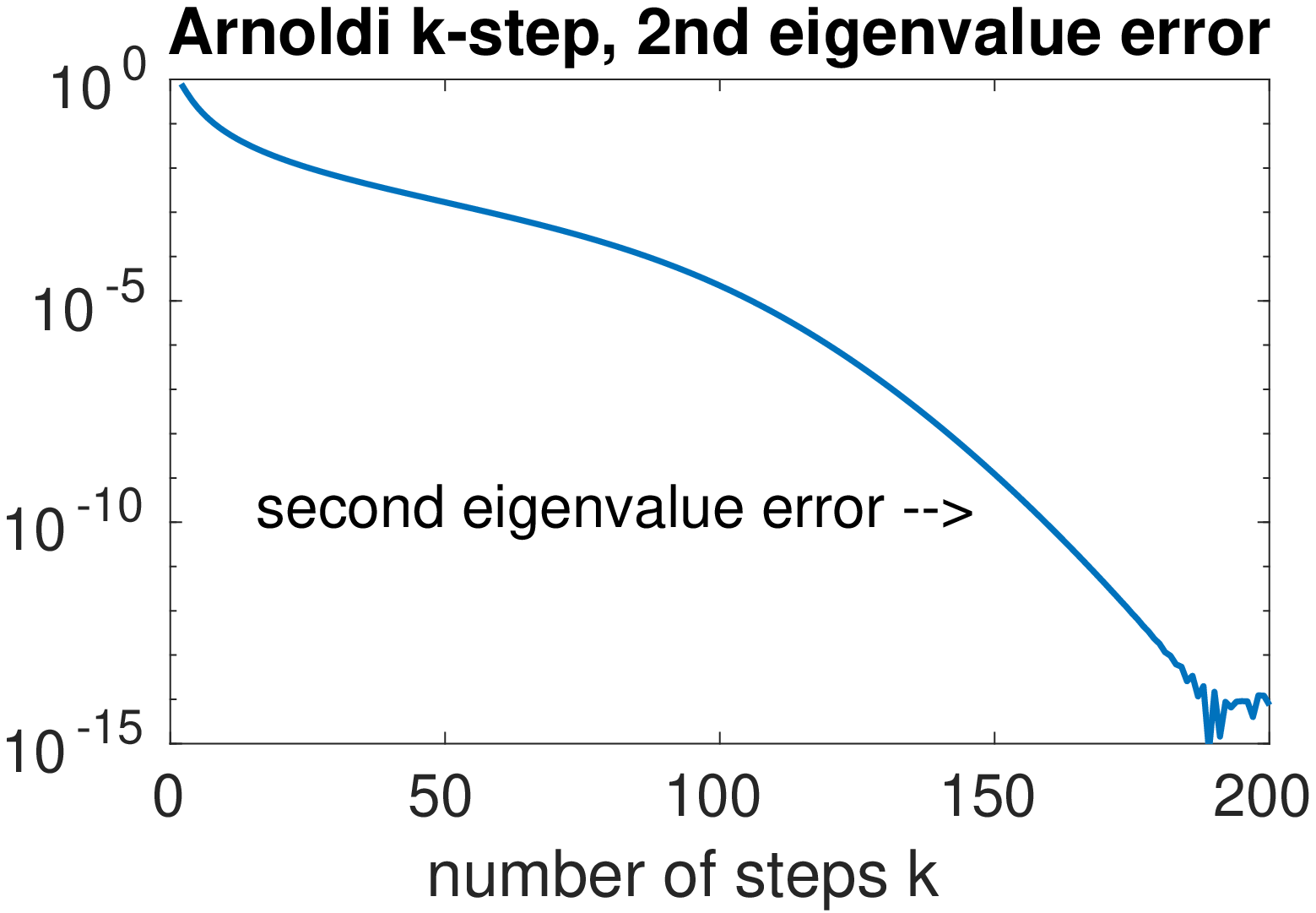}}
\caption{The Arnoldi $k$-step method for the $n\times n$ diagonal matrix $A$
where $a_{ii}=n/i$, $i=1,\dots,n$, for $n=1000$.
(a) Errors for the largest eigenvalue (which is 1) and corresponding eigenvector
residual norm.
(b) Errors for the second-largest eigenvalue (which is $1-1/n$).}
\label{fig:plotkstep}
\end{figure}

One can show by induction that the vectors $\yy_j$ are orthogonal, so that,
in exact arithmetic, $H=K$.
But this subtle change makes the algorithm far more robust, as shown in 
Figure \ref{fig:plotkstep}(a).
The $k$-step algorithm approximates accurately many of the largest eigenvalues.
Figure \ref{fig:plotkstep}(b) plots the error for the second-largest eigenvalue.
Figure \ref{fig:plotkstep}(a) further displays the expected \cite{lrsBIBhi} behavior 
that the eigenvalue error is proportional to the square of the eigenvector error. 

\subsection{Restarted Arnoldi algorithm}

Our $k$-step algorithm described above
is a restarted Arnoldi algorithm \cite{ref:sorensensurvey}.
What is different is that here we propose the restart for purely algorithmic
purposes (e.g., minimizing storage) as opposed to deflation 
or stability considerations \cite{ref:imprestartArnoldi}.

\subsection{Rationale for $k$-step algorithm}

Assessments of efficiency are problem dependent, and so we look in detail
at one application.
Suppose that the matrix $A$ is sparse and of size $N\times N$.
We quantify the sparseness by assuming that a matrix-vector multiplication
$A\yy$ costs $N\varkappa$ operations.
Vector dot products cost $N$ operations.
Thus the Arnoldi method requires 
$$
kN\varkappa + k^2 N
$$
operations.
In many cases, multiplication by $A$ actually requires solving a system of
equations \cite{lrsBIBis}, and so $\varkappa$ might be as large as 1000 or larger, 
and perhaps increasing as $N$ increases for solvers that are not optimal order.
And in many physical applications, only a few eigenvalues are of interest.

One limit on $k$ is the required $kN$ storage of the Krylov vectors.
Thus there may be a need for several iterations of a $k$-step method for a fixed $k$, 
even though one iteration for a much larger $k$ might be more efficient.

The restarted algorithm requires solving a $k \times k$ eigenproblem at
each step.
Let us assume that the reduced eigenproblem takes on the order of $k^3$ operations.
Then if $k<<\min\{\varkappa,N\}$, the cost of the reduced eigensolve is negligible.

\section{Accelerating the $k$-step method with extrapolation}

A depth-1 extrapolation can be applied to the $k$-step or restarted Arnoldi
method as a simple low-cost post-processing procedure.
As discussed and shown below, this can be
advantageous particularly for smaller values of $k$.
Here we discuss the extrapolation with iteration $j$ parameter $\gamma_j$ as summarized 
in the following algorithm. In order for the first two approximate eigenvalues
$\lambda_1^{(j)}$ and $\lambda_2^{(j)}$ to provide meaningful information,
generally $k \ge 4$ makes sense.  However, if the second approximate eigenvalue
is not being used to define the extrapolation parameter $\gamma_j$, then $k\ge 2$
makes sense.

\begin{alg}[Extrapolated $k$-step Arnoldi]\label{alg:kstepex} 

\noindent Choose $y^{(0)}$ and  $k\ge2$.
\\ \noindent
Compute $[y^{(1)},\lambda_1^{(0)},\lambda_2^{(0)}]$ 
= Arnoldi$(y^{(0)},A,k)$. Set $u^{(1)} = y^{(1)}$.
\\ \noindent
for $j = 1, 2, \ldots$
\\ \indent
a. Compute $[y^{(j+1)},\lambda_1^{(j)},\lambda_2^{(j)}]$ = Arnoldi$(u^{(j)},A,k)$
\\ \indent
b. Set $\gamma_j$
\\ \indent
c. Set $u^{(j+1)} = (1 - \gamma_j) y^{(j+1)} + \gamma_j y^{(j)}$
\\ \noindent
end
\end{alg}
In the tests of subsections \ref{subsec:tests} and \ref{subsec:details}, the 
condition to exit the loop on convergence is 
$\norm{A y^{(j+1)} - \lambda_1^{(j)} y^{(j+1)}} <$ {\tt tol}, for a given tolerance
{\tt tol}.

The choice of extrapolation parameter $\gamma_j$ is the key to a successful 
extrapolation. 
In \cite{ref:Pollackextraeig}, where the power iteration is accelerated by 
extrapolation, the parameter $\gamma_j$ which gives asymptotically exponential
convergence for positive semidefinite problems is an approximation of 
$-(\lambda_2/\lambda_1)^j$, where the eigenvalues of $A$ are labeled in descending 
magnitude.  Here we will see the approximation of $-(\lambda_2/\lambda_1)^j$
produced by the $k$-step Arnoldi method
gives an effective acceleration of the $k$-step method; but, it is not clear that this
is necessarily the best choice.

The complication in setting the extrapolation parameter
lies in understanding the approximate eigenvector $y^{(j+1)}$ produced
from Arnoldi$(u^{(j)},A,k)$ as an expansion in the eigenbasis of $A$; in contrast to
the power iteration, this expansion is not available in closed form for the $k$-step
method for general values of $k$.
For concreteness, suppose $A$ is diagonalizable with a basis of orthonormal eigenvectors
$\{v_i\}_{i = 1}^n$, corresponding to eigenvalues $\{\lambda_i\}_{i = 1}^n$, 
labeled with decreasing magnitude.

\subsection{Analysis of extrapolation}

By construction, the first generated iterate 
$y^{(1)} \in \cK_k(y^{(0)}) = \spa\{y^{(0)}, Ay^{(0)}, \ldots A^{k-1}y^{(0)}\}$, the
$k$-dimensional Krylov space generated by $A$ applied to $y^{(0)}$.
The next iterate $y^{(2)} \in \cK_k(y^{(1)}) \subsetneq \cK_{2k-1}(y^{(0)})$, and 
in general $y^{j+1} \in \cK_k(y^{(k)}) \subsetneq \cK_{(j+1)k -j}(y^{(0)})$.

A formal expansion of each $y^{(j)}$ in terms of the eigenbasis of $A$ can be 
expressed as
\[
y^{(j)} = \frac{1}{h_j} \sum_{i = 1}^n p_i^{(j)}(\lambda_i) v_i,
\quad h_j = \bigg( \sum_{i = 1}^n (p_i^{(j)}(\lambda_i))^2 \bigg)^{1/2},
\]
where $p_i^{(j)}(\lambda_i)$ is a polynomial of degree at most $kj-(j-1)$ in $\lambda_i$.
Now let's consider the ratio of the components of $u^{(j+1)}$ and $y^{(j)}$ 
in the direction of each eigenvector $v_i$.
First define 
\begin{align}\label{eqn:effi-a}
\eta_i^{(j+1)} \coloneqq
\frac{v_i \cdot \proj_{v_i} y^{(j+1)}}
     {v_i \cdot \proj_{v_i} y^{(j)}  } =
\frac{h_{j}}{h_{j+1}} \frac{p_i^{(j+1)}(\lambda_j)}{p_i^{(j)}(\lambda_j) } 
= 
\frac{p_i^{(j+1)}(\lambda_j)/h_{j+1}}{p_i^{(j)}(\lambda_j)/h_j }. 
\end{align}
Then, noting 
\[
\proj_{v_i}u^{(j+1)} = \frac{p_i^{(j)}(\lambda_j)}{h_j}v_i
\Big( (1-\gamma_j) \frac{h_j}{h_{j+1}} 
\frac{p_i^{(j+1)}(\lambda_i)}{p_i^{(j)}(\lambda_i)} + \gamma_j
\Big),
\]
we have
\begin{align}\label{eqn:effi}
\widehat \eta_i^{(j+1)}\coloneqq
\frac{v_i \cdot \proj_{v_i} u^{(j+1)} }
     {v_i \cdot \proj_{v_i} y^{(j)}   } =
(1 - \gamma_j)\eta_i^{(j+1)} + \gamma_j.
\end{align}

\begin{table}
\centering
{\renewcommand{\arraystretch}{1.1}
\begin{tabular}{|l||c|c|c|c|c|c|c|}
\hline
\backslashbox{Matrix}{$\gamma_j$}
& 0 & $-0.25$ & $-0.5$ & $-0.75$ & 
          $ -\left|\frac{\lambda_2^{(j)}}{2\lambda_1^{(j)}}\right|^2$ &
          $ -\left|\frac{\lambda_2^{(j)}}{\lambda_1^{(j)}}\right|$ &
          $ -\left|\frac{\lambda_2^{(j)}}{\lambda_1^{(j)}}\right|^j$ 
\\ [2pt]
\hline
$A_1$                    & 192 & 94 & 73 & 76 & 80 & 97 & 98 \\
$A_2 = ${\tt Kuu}        & 86  & 47 & 53 & 54 & 43 & 56 & 46 \\
$A_3 = ${\tt ifiss\_mat} & 165 & 105& 83 & 42 & 79 & 72 & 68 \\
$A_4 = ${\tt gearbox}    & 157 & 48 & 52 & 52 & 56 & 82 & 82 \\
$A_5 = ${\tt ss1}        & 85  & 93 & 91 & 95 & 212& 390& 75 \\
$A_6 = ${\tt Si87H76}    & 63  & 32 & 37 & 37 & 32 & 37 & 33 \\

\hline
\end{tabular}
}
\vspace{4mm}
\caption{Number of iterations to residual convergence of $10^{-7}$ for constant and
dynamically chosen extrapolation parameters used in Algorithm \ref{alg:kstepex} 
with $k=8$.}
\label{tab:extrap-param}
\end{table}

Now we are interested in how the modes defined by \eqref{eqn:effi} grow or decay 
compared to those defined by \eqref{eqn:effi-a}. 
Restricting $\gamma_j$ to the interval
$[-1,0]$ damps the modes for which $\eta_i^{(j+1)}>0$.
Notice that if $\gamma_j = 1$, then $\widehat \eta_i^{(j+1)} = \gamma_j = 1$, 
yielding complete stagnation: the last iterate was repeated.  
If $\gamma_j = 0$, then $\widehat \eta_i^{(j+1)} = \eta_i^{(j+1)}$ for all
$i$ (no extrapolation was applied).  
If $\gamma_j = -\eta_i^{(j+1)}$, then 
$\widehat \eta_i^{(j+1)} = (\eta_i^{(j+1)})^2$ which quadratically accelerates decay 
of mode $i$;
 and, if $\gamma_j = -1$, then $\widehat \eta_i^{(j+1)} = 2 \eta_i^{(j+1)} -1$, 
which causes rapid decay of modes for which $\eta_i^{(j+1)} \approx 0.5$, and slow
decay of modes for which $\eta_i^{(j+1)}$ is close to zero or one.
More generally, for $\gamma_j < 0$ it holds that 
$\widehat \eta_i^{(j+1)} <  \eta_i^{(j+1)}$ so long as $\eta_i^{(j+1)} > 0$;
and, $\widehat \eta_i^{(j+1)} > -\eta_i^{(j+1)}$ (no growth in magnitude) for
$\eta_i^{(j+1)} < \gamma_j/(\gamma_j -2)$.

We are most interested in the case of $\eta_i^{(j+1)} > 0$, as both $p_i^{(j+1)}$ and
$p_i^{(j)}$ will tend to the same sign, particularly for the $k/2$ largest modes
(assuming $k$ even), as these are the eigenvalues better approximated using the
$k$-step method.  However, changes in sign of $p_i^{(j)}$ can occur between iterations,
particularly in the $k/2$ smallest modes.
For $\eta_i^{(j+1)} < 0$, letting $\gamma_j < 0$ always causes
$\widehat \eta_i^{(j+1)} < \eta_i^{(j+1)}$. Growth in mode $i$ occurs 
if $\eta_i^{(j+1)}< -(1 + \gamma_j)/(1-\gamma_j)$: increased growth for
$-1 <\gamma_j < 0$ larger in magnitude. Since the Arnoldi or $k$-step process in 
itself is efficient at quickly damping the smaller modes, the purpose of the
extrapolation is to accelerate the decay of the larger subdominant modes that 
compete with the dominant mode, hence we will consider extrapolation parameters
restricted to the interval $[-1,0]$.

\subsection{Convergence tests}\label{subsec:tests}

\begin{figure}
\centering
\includegraphics[trim = 0pt 0pt 10pt 0pt,clip = true, width=0.32\textwidth]
{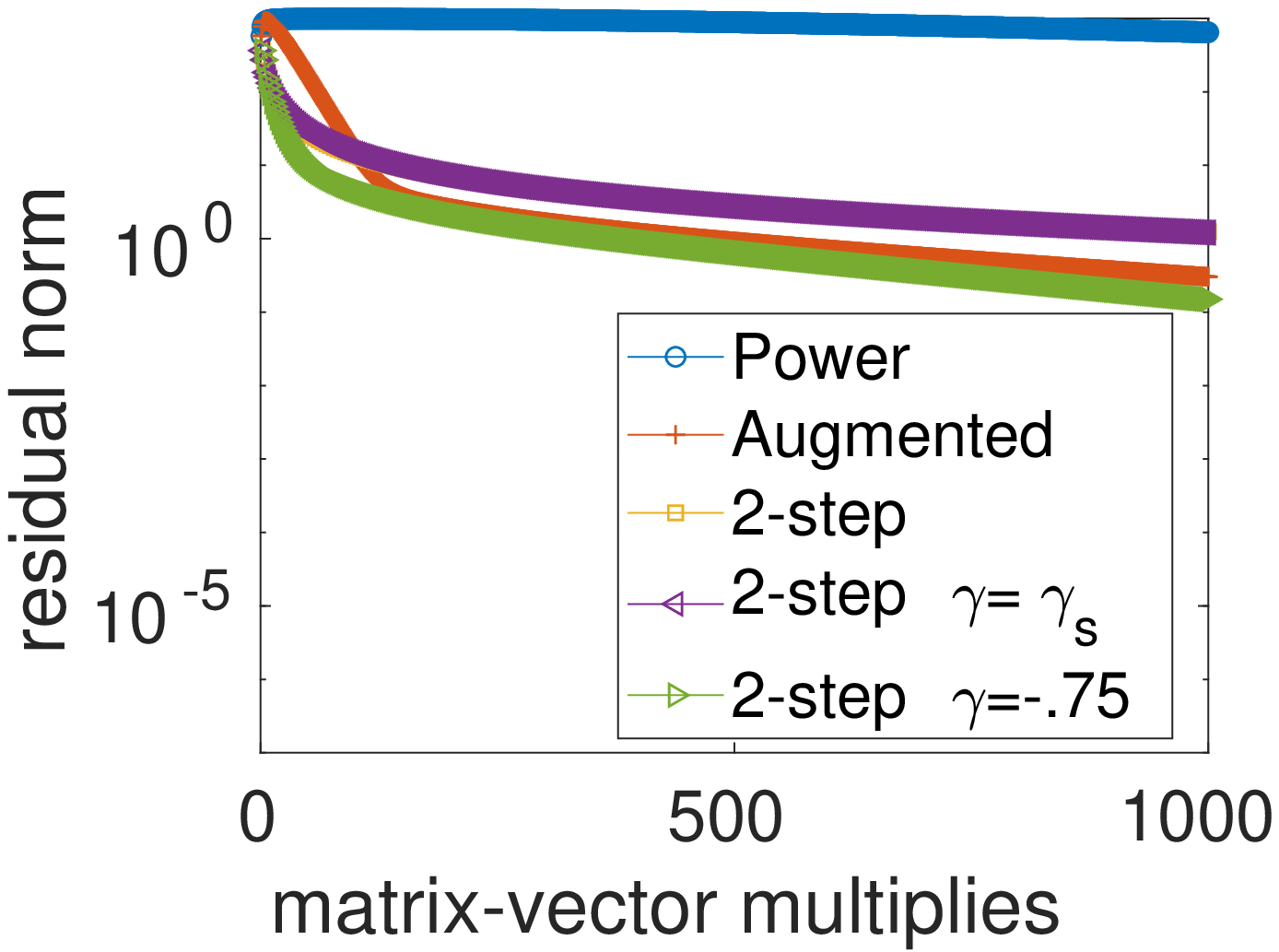}
\includegraphics[trim = 0pt 0pt 10pt 0pt,clip = true, width=0.32\textwidth]
{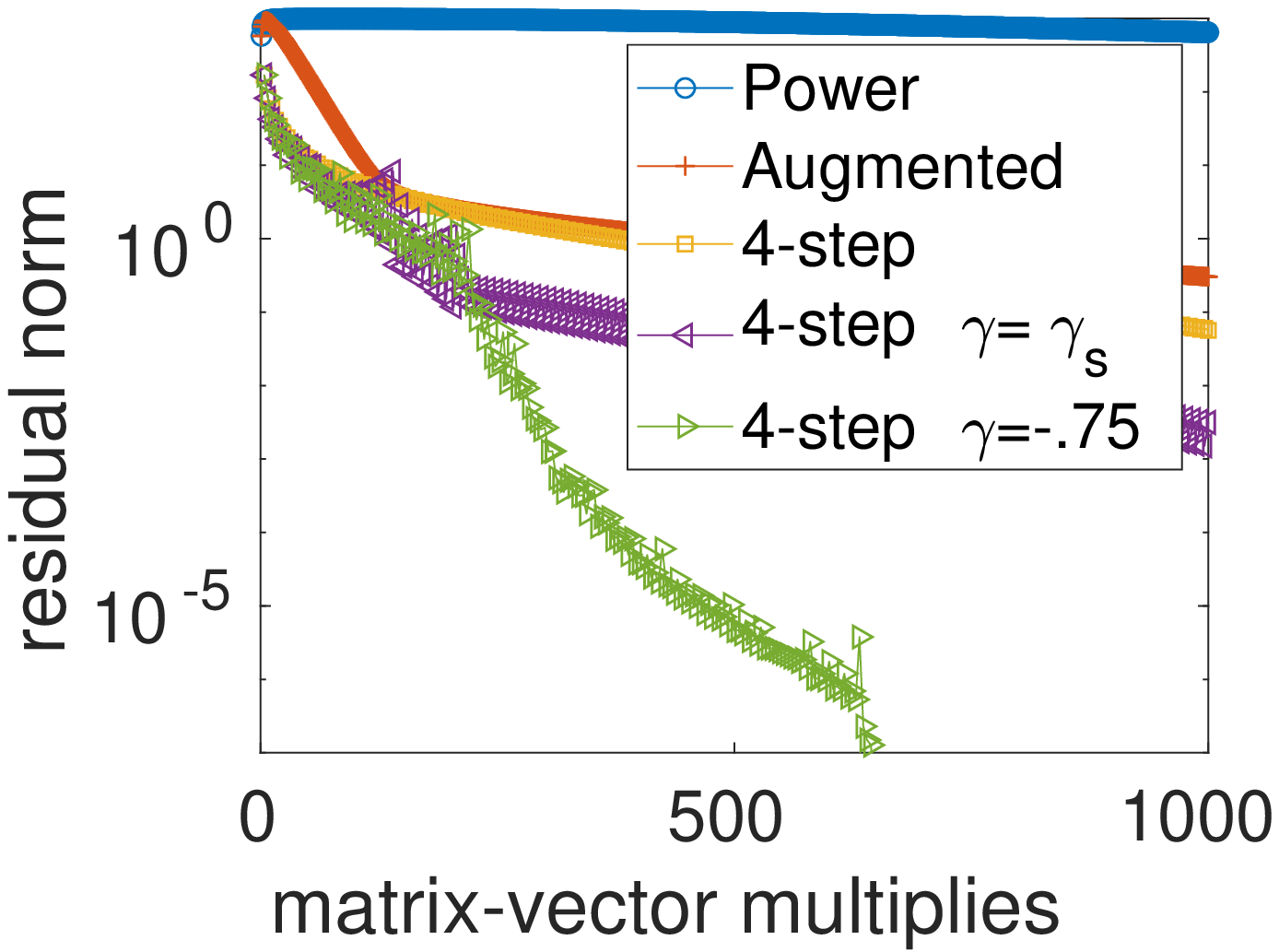}
\includegraphics[trim = 0pt 0pt 10pt 0pt,clip = true, width=0.32\textwidth]
{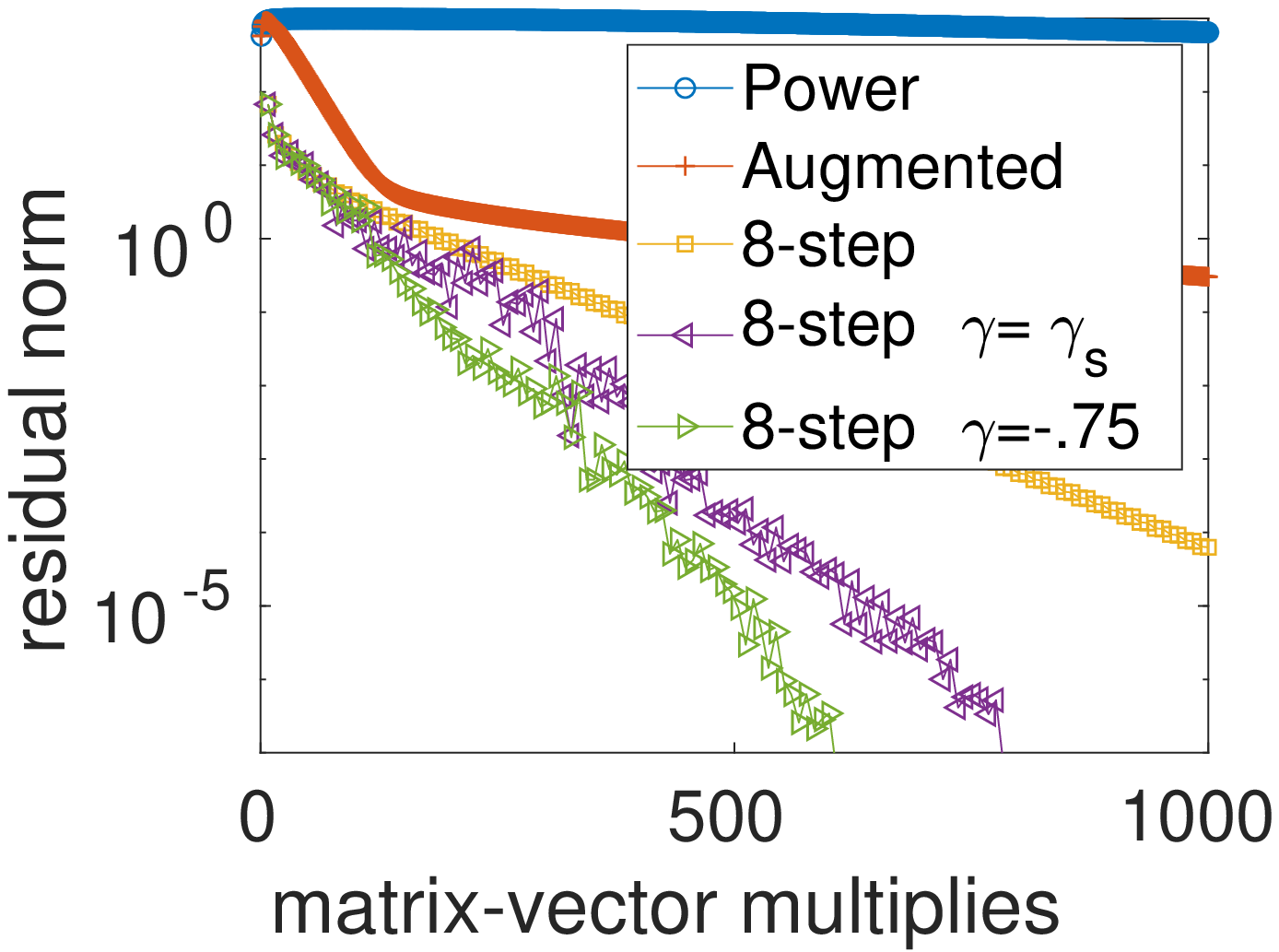}
\caption{Residual history for Algorithm \ref{alg:kstepex}, with $k = 2$ (left), 
$k=4$ (center) and $k=8$ (right), with $\gamma_j = 0$ ($k$-step),
$\gamma_j = \gamma_s =  -\big(\lambda_2^{(j)}/\lambda_1^{(j)}\big)^j$
and $\gamma_j = -0.75$, compared with the power iteration and Augmented
method of \cite{ref:Pollackextraeig} with $\eta=40$,
applied to $A_1 = \diag(1000,-999,\ldots, 2,-1)$.
}
\label{fig:idiag-extrap}
\end{figure}

In Table \ref{tab:extrap-param}, the number of iterations to convergence with six 
choices of extrapolation parameter $\gamma_j$ are compared on the following test set:
$A_1 = \diag(1000,-999,998,\ldots,2,-1)$;  $A_2 = ${\tt Kuu}, symmetric positive definite
with $n = 7102$; $A_3 = $ {\tt ifiss\_mat}, neither positive definite nor
symmetric, with $n = 96307$; $A_4 = $ {\tt gearbox}, symmetric and not positive definite,
with $n = 153746$; $A_5 =$ {\tt ss1}, neither positive definite nor symmetric, with
$n = 205282$; $A_6 = $ {\tt Si87H76}, symmetric and not positive definite, with
$n = 240369$.
Matrices $A_2$---$A_6$ are available from the
SuiteSparse matrix collection \cite{DH11}. Each iteration is started with 
$y^{(0)} = (1, \ldots, 1)^t$. 

Three constant extrapolation parameters, $\gamma_j =-0.25,-0.5,-0.75$ are compared
with three dynamically chosen parameters, 
$\gamma_j = -|\lambda_2^{(j)}/\lambda_1^{(j)}|^2/4$,
$\gamma_j = -|\lambda_2^{(j)}/\lambda_1^{(j)}|$ and 
$\gamma_j = -|\lambda_2^{(j)}/\lambda_1^{(j)}|^j$,
where $\lambda_1^{(j)}$ and $\lambda_2^{(j)}$ 
are the first and second eigenvalues returned on iteration $j$ from
Arnoldi$(u^{(j)},A,k)$.  
The parameter $\gamma_j = -|\lambda_2^{(j)}/\lambda_1^{(j)}|^2/4$ is related to that
used in \cite{DSHMRX19}.
The parameters 
$\gamma_j = -|\lambda_2^{(j)}/\lambda_1^{(j)}|$ and
$\gamma_j = -|\lambda_2^{(j)}/\lambda_1^{(j)}|^j$,
are based on the extrapolation parameter used in \cite{ref:Pollackextraeig}.
The last of these,
$\gamma_j = \gamma_s = -|\lambda_2^{(j)}/\lambda_1^{(j)}|^j$
is denoted with subscript ``$s$'' due to its similarity to the parameter that defines
the ``simple'' method in \cite{ref:Pollackextraeig}, which 
has established convergence properties for power iterations.
All three dynamically chosen parameters are ensured to satisfy $-1 \le \gamma_j \le 0$, 
which is the interval of interest, as per the discussion above.
On the whole, the constant extrapolation parameters
significantly reduce computation in most cases, with a minor increase in the number
of iterations in one tested case. The dynamically chosen parameter
$\gamma_j = -|\lambda_2^{(j)}/\lambda_1^{(j)}|^j$ reduced the number of iterations in
all cases, although not as effectively as the constant extrapolation in the case
of $A_4 = ${\tt gearbox}.
With the exception of $A_5 = ${\tt ss1}, the extrapolation methods reduced the
number of iterations (without any additional matrix-vector multiplies)
generally by 50\% or more. 
The two extrapolation parameters from table \ref{tab:extrap-param} which provide
the greatest benefit in most cases, while causing minimal damage in the worst, 
are the constant parameter
$\gamma_j = -0.75$, and the dynamic paramter 
$\gamma_j = -|\lambda_2^{(j)}/\lambda_1^{(j)}|^j$.  These two successful 
parameters are next considered for different choices $k$ in the $k$-step method.

\begin{figure}
\centering
\includegraphics[trim = 0pt 0pt 10pt 0pt,clip = true, width=0.32\textwidth]
{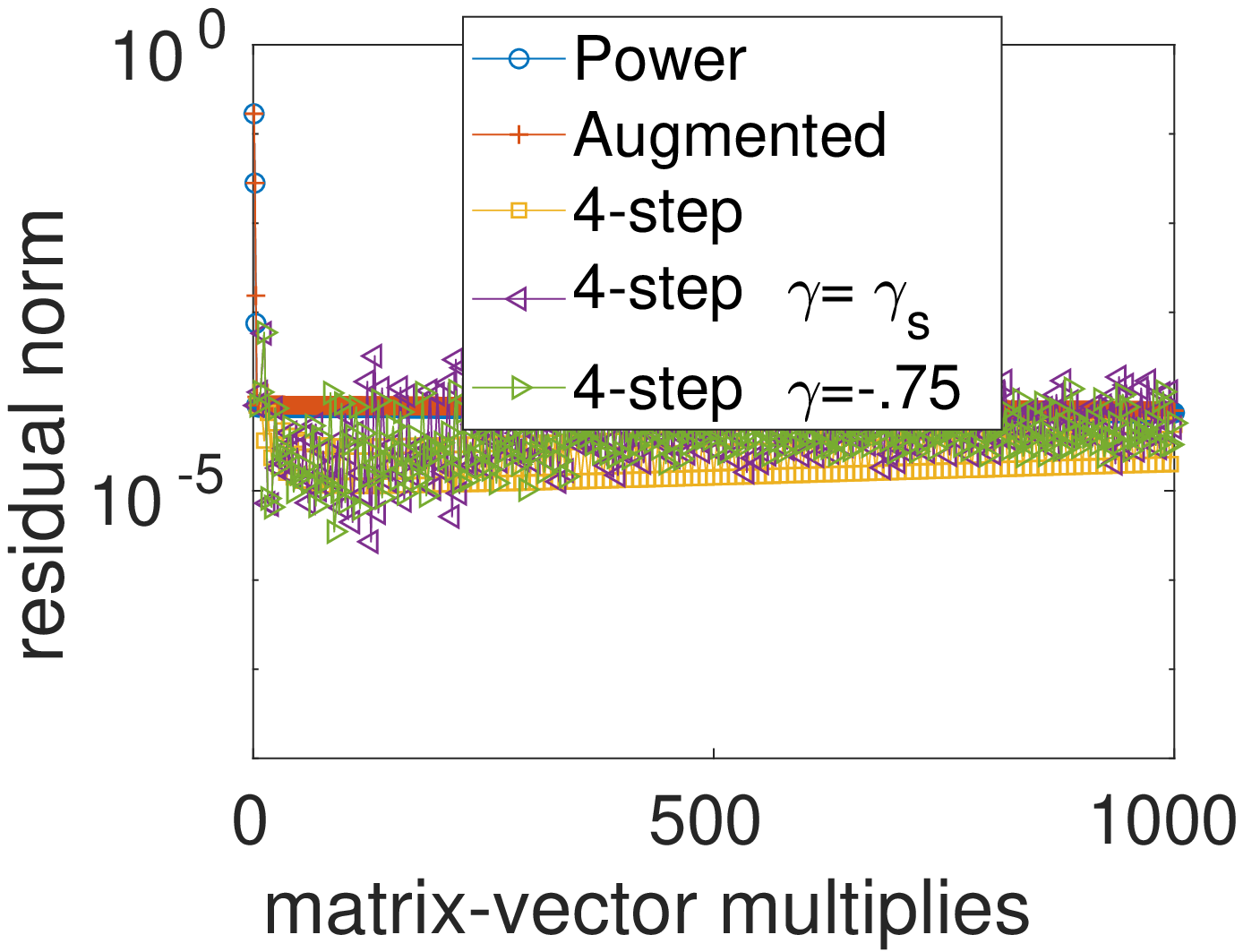}
\includegraphics[trim = 0pt 0pt 10pt 0pt,clip = true, width=0.32\textwidth]
{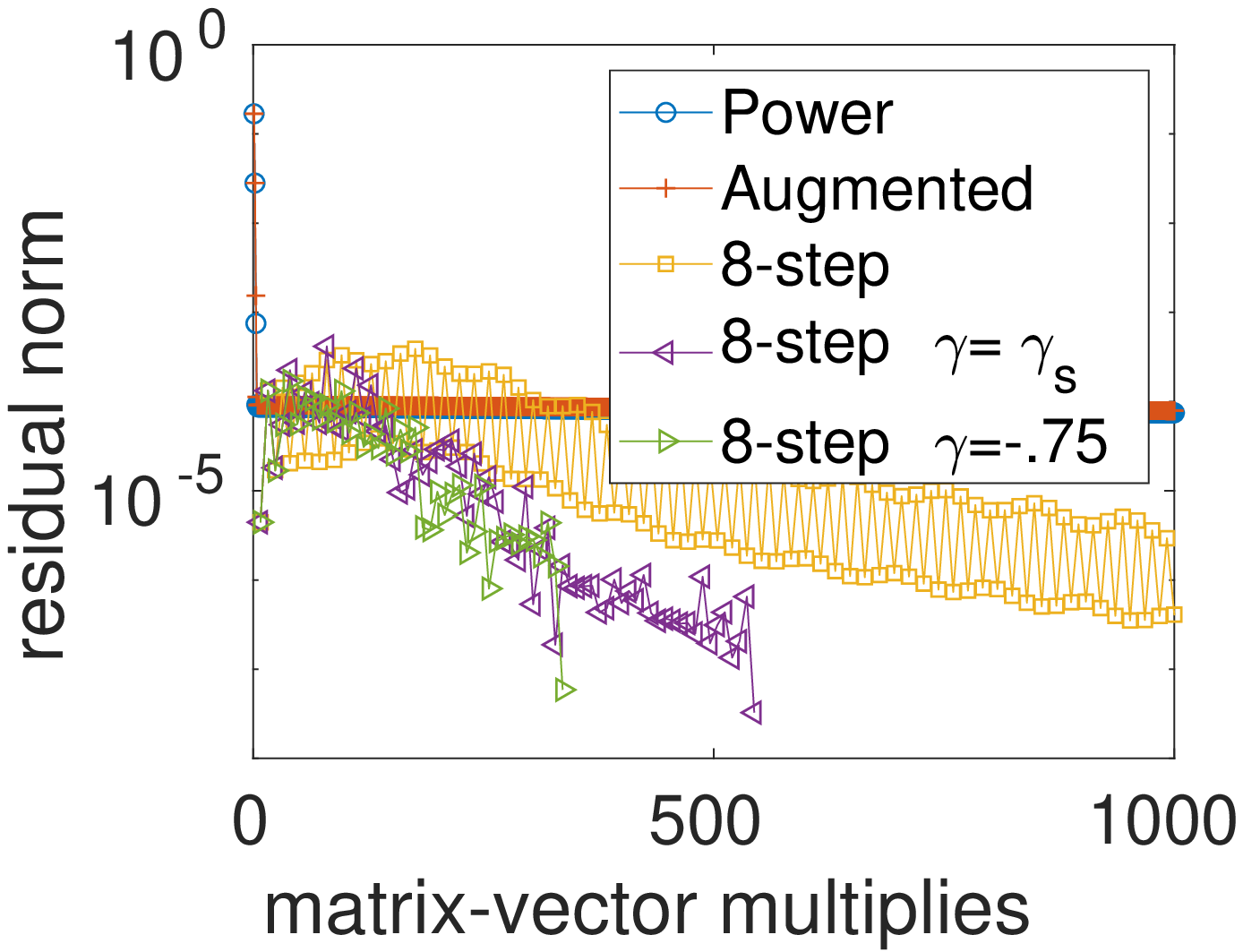}
\includegraphics[trim = 0pt 0pt 10pt 0pt,clip = true, width=0.32\textwidth]
{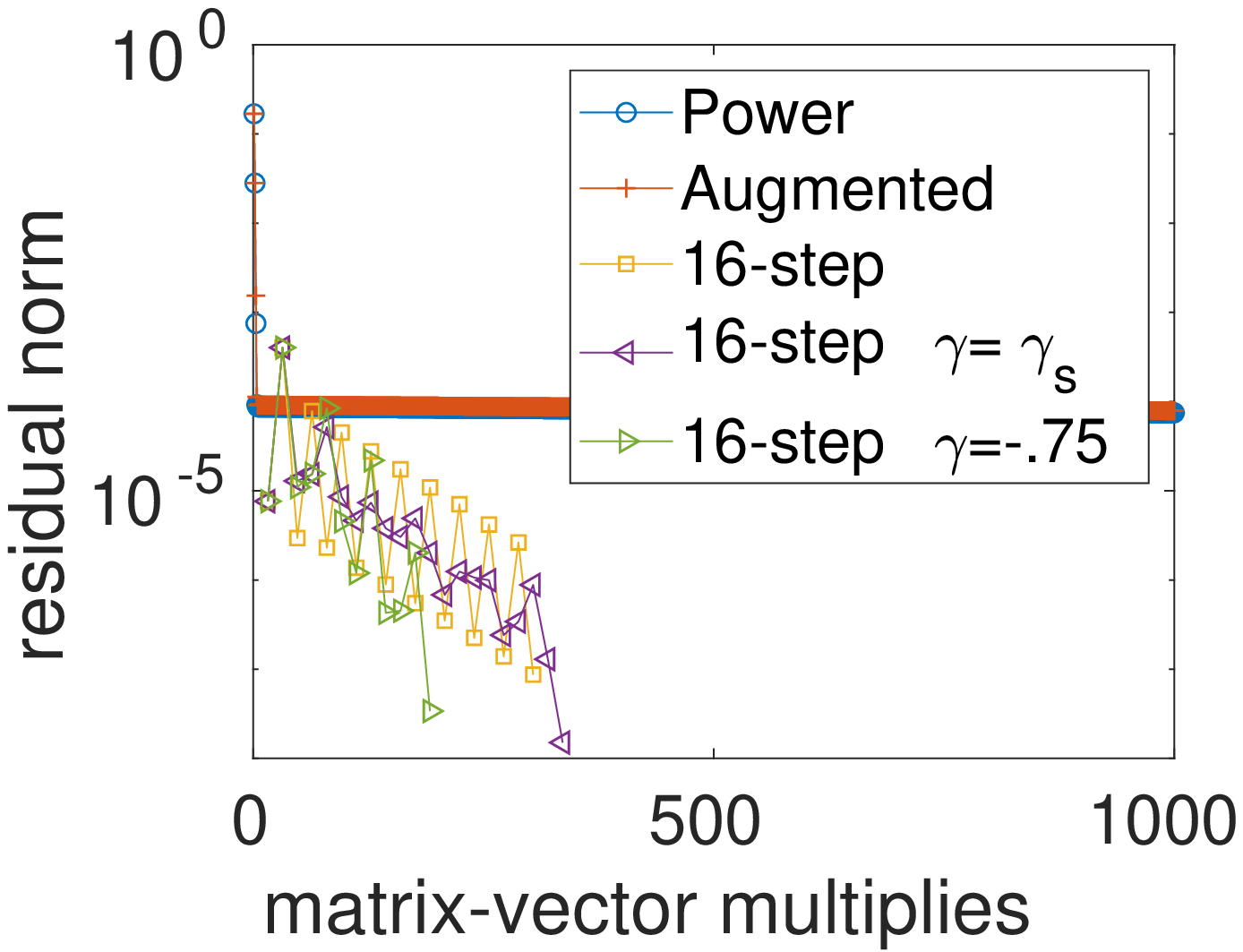}
\caption{Residual history for Algorithm \ref{alg:kstepex}, with $k = 4$ (left), 
$k=8$ (center) and $k=16$ (right), with $\gamma_j = 0$ ($k$-step), 
$\gamma_j = \gamma_s = -\big(\lambda_2^{(j)}/\lambda_1^{(j)}\big)^j$
and $\gamma_j = -0.75$, compared with the power iteration and Augmented
method of \cite{ref:Pollackextraeig} with $\eta=40$, applied to $A_3 = $ {\tt ifiss\_mat}.
}
\label{fig:ifiss-extrap}
\end{figure}

\subsection{Convergence details}\label{subsec:details}

\begin{figure}
\centering
\includegraphics[trim = 0pt 0pt 10pt 0pt,clip = true, width=0.32\textwidth]
{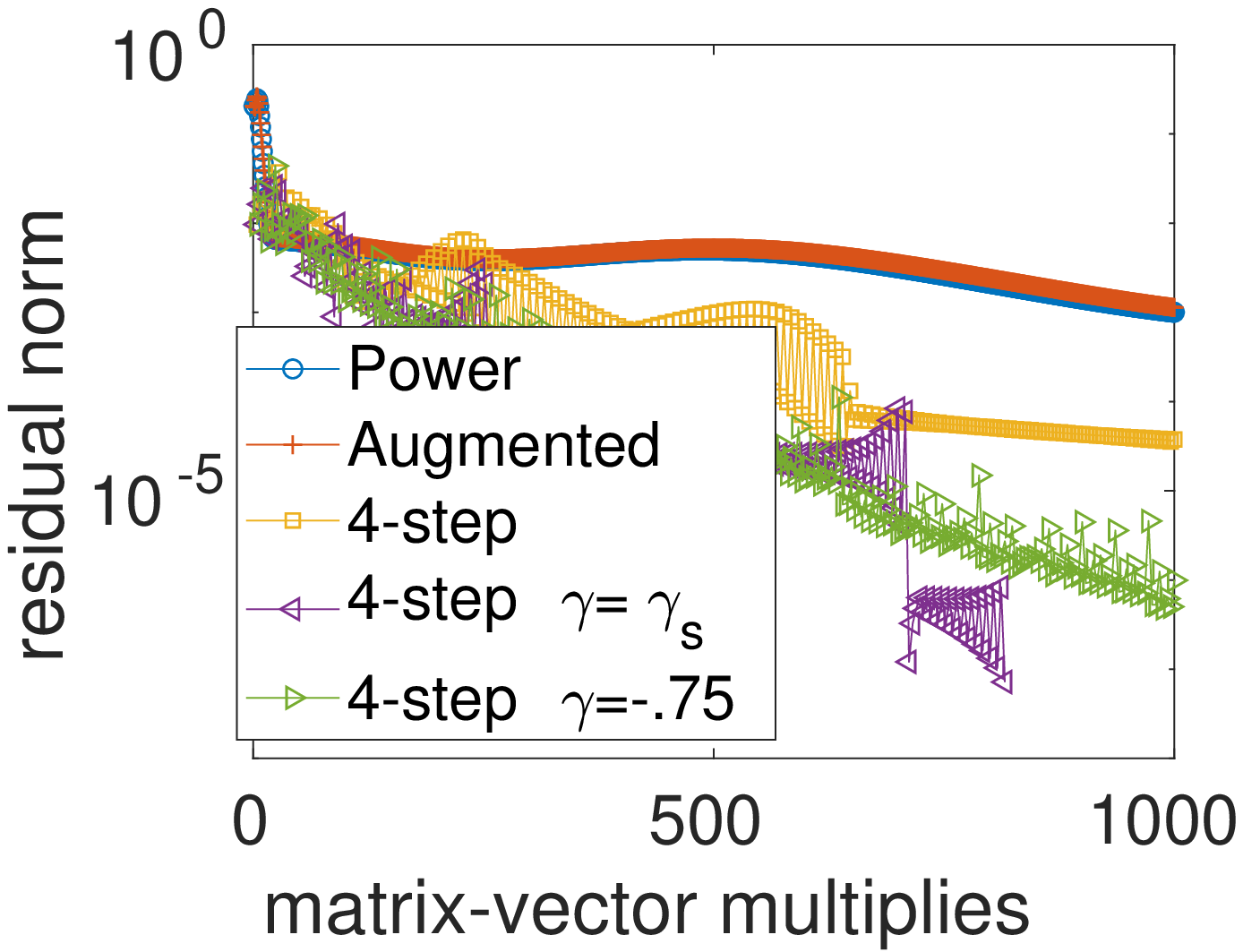}
\includegraphics[trim = 0pt 0pt 10pt 0pt,clip = true, width=0.32\textwidth]
{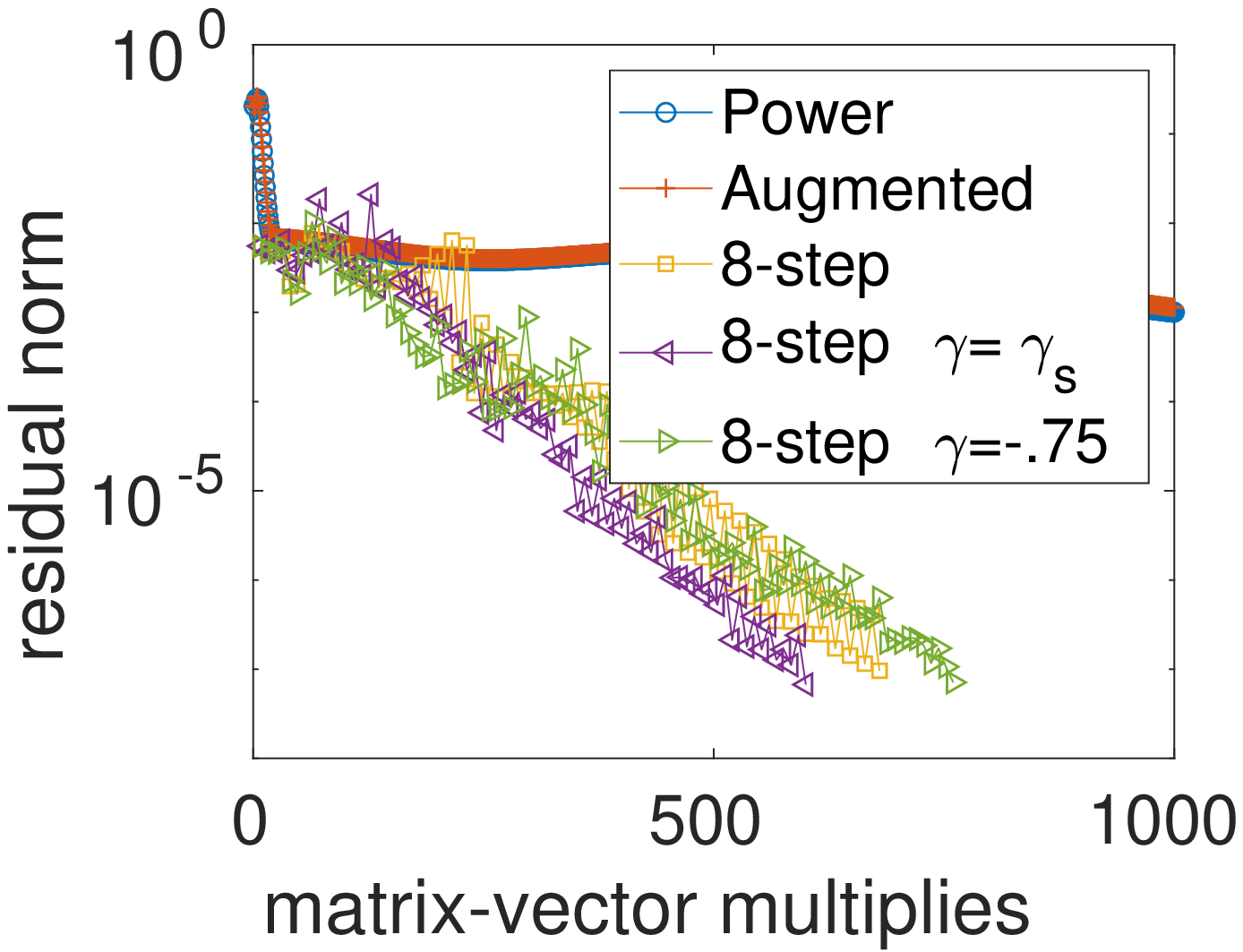}
\includegraphics[trim = 0pt 0pt 10pt 0pt,clip = true, width=0.32\textwidth]
{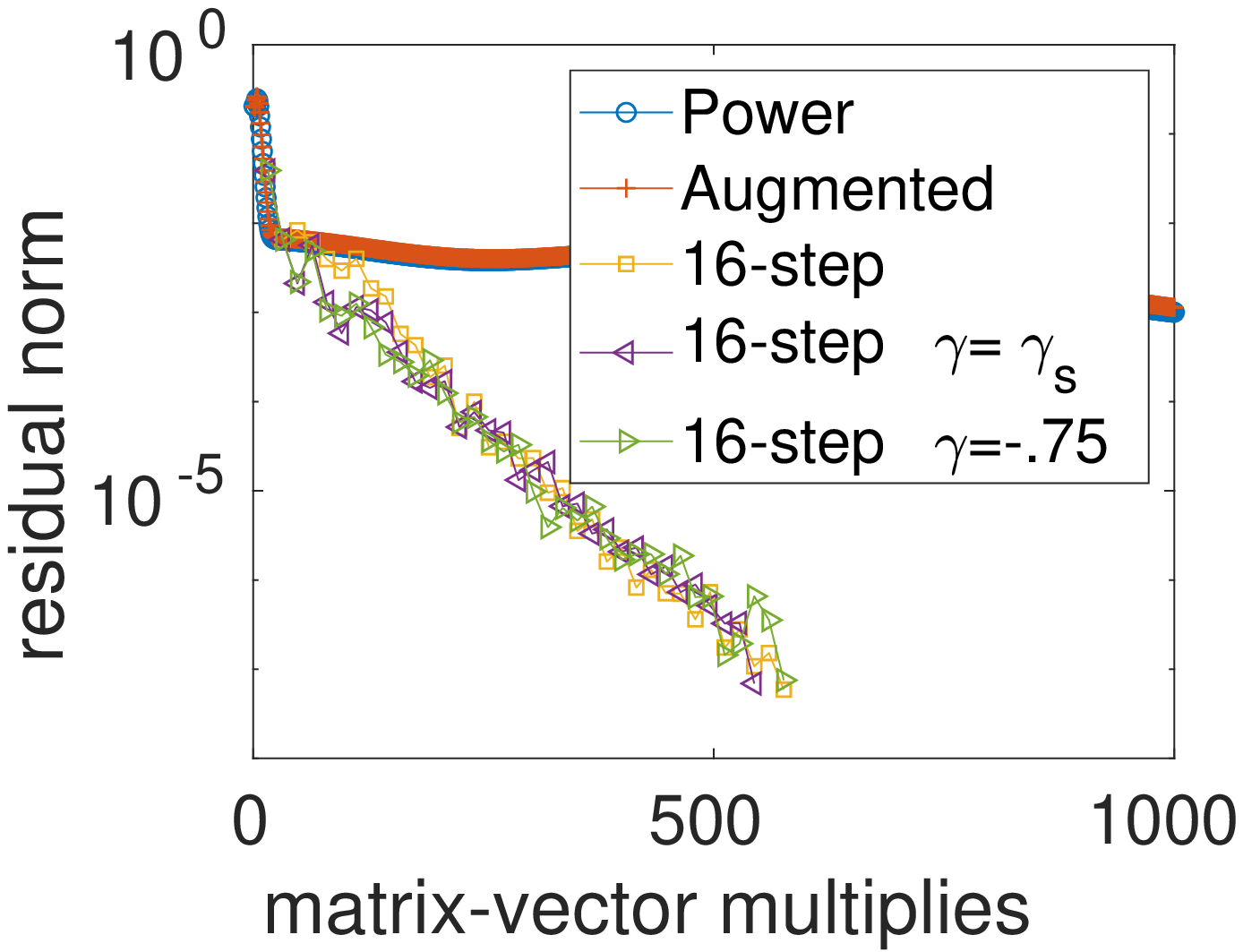}
\caption{Residual history for Algorithm \ref{alg:kstepex}, with $k = 2$ (left),
$k=4$ (center) and $k=8$ (right), with $\gamma_j = 0$ ($k$-step),
$\gamma_j = \gamma_s = -\big(\lambda_2^{(j)}/\lambda_1^{(j)}\big)^j$
and $\gamma_j = -0.75$, compared with the power iteration and Augmented
method of \cite{ref:Pollackextraeig} with $\eta=40$, applied to $A_5 =$ {\tt ss1}.
}
\label{fig:ss1-extrap}
\end{figure}

Figures \ref{fig:idiag-extrap}---\ref{fig:ss1-extrap} compare 
Algorithm \ref{alg:kstepex} with parameters $\gamma_j = -0.75$ and 
$\gamma_j =\gamma_s =  -|\lambda_2^{(j)}/\lambda_1^{(j)}|^j$ with the $k$-step method
($\gamma_j = 0$) using
different values of $k$, together with
the power iteration and the Augmented mode-damping method of \cite{ref:Pollackextraeig}.
Each plot shows the $l_2$ norm of the residual for the dominant eigenpair on the 
$y$-axis, and the number of matrix-vector multiplies on the $x$-axis.
In figure \ref{fig:idiag-extrap} it is notable for matrix $A_1$
that the constant extrapolation
$\gamma_j = -0.75$ accelerates convergence for $k = 4,8$, whereas the dynamically
assigned parameter $\gamma_s$ accelerates convergence nearly as well, but only for
$k=8$. This illustrates that the constant parameter may be preferable in cases where
the second eigenvalue may not be approximated well. Figure \ref{fig:ifiss-extrap}
(left) with $k=4$ illustrates with $A_3 = $ {\tt ifiss\_mat}
that the extrapolation does not necessarily induce
convergence where the $k$-step method itself does not converge. The center and
right plots with $k=8$ and $k=16$ illustrate that the extrapolation can still
significantly improve convergence where the residual for the $k$-step method
demonstrates oscillatory behavior.
Figure \ref{fig:ss1-extrap} (left) with $k=4$ illustrates
with $A_5 = $ {\tt ss1} an instance where the
dynamically chosen $\gamma_s$ reduces the number of iterations more effectively than
the constant parameter for a smaller value of $k$, however the center and right plots
with $k = 8$ and $k=16$ show little difference with or without the extrapolation.

The examples shown in table \ref{tab:extrap-param} and figures
\ref{fig:idiag-extrap}---\ref{fig:ss1-extrap} together illustrate that for large enough
values of $k$, the $k$-step
method converges efficiently with respect to the number of matrix-vector multiplies;
however, for the range of $k$ values for which the $k$-step method converges but not
efficiently, the extrapolation improves convergence and restores efficiency.
These results are consistent with the observations in \cite{lrsBIBjf}
and \cite{saad1980arnoldi} that the
restarted Arnoldi or $k$-step method with $k = j\times k'$ is more effective than $j$
iterations of the $k'$-step method.  However, the extrapolation reduces the
sensitivity of the method to the choice of $k$, so that smaller values of $k$ can produce
nearly the same efficiency, and with less demand on system memory.

\section{LOBPCG and $k$-step methods}

Locally Optimal Block Preconditioned Conjugate Gradient (LOBPCG) methods
\cite{benner2011locally,knyazev2017recent,knyazev2001toward}
have proved to be very effective for solving eigenproblems.
Here we compare them with the $k$-step methods.
For simplicity, we restrict to the case that there is no preconditioning
and the eigenproblem is standard, not generalized.
Referring to the Wikipedia page\footnote{{\tt https://en.wikipedia.org\_wiki\_LOBPCG}}
for LOBPCG, we can write the basic algorithm as follows.
First define the Rayleigh quotient $\rho(\xx)$ by
$$
\rho(\xx)=\frac{\xx^tA\xx}{\xx^t\xx},
$$
and the residual
$$
\rr(\xx)=A\xx-\rho(\xx)\xx.
$$
Given $\xx_i$, LOBPCG defines
$$
\xx_{i+1}=\arg{} \max_{\yy\in\,\hbox{\small span}\{\xx_i,\rr(\xx_i)\}} \rho(\yy).
$$
But assuming that $A\xx_i\neq\rho(\xx_i)\xx_i$, that is, the algorithm has not converged,
$$
\hbox{span}\{\xx_i,\rr(\xx_i)\} =\hbox{span}\{\xx_i,A\xx_i \}.
$$
Therefore this form of LOBPCG is the same as the 2-step algorithm.
In the Wikipedia page for LOBPCG, the 3-step algorithm is also described.

\section{Enhancements of $k$-step methods}

There are several possible enhancements of $k$-step methods.

\subsection{Block $k$-step methods}

Block $k$-step methods can also be useful.
For simplicity, suppose that $k$ is even.
A block 2-step method, with block size $k/2$, would work as follows.
First pick $k/2$ initial vectors $\ww_i$.
(This could be done by generating $k/2$ Krylov vectors.)
Next perform one Krylov step for each $i$:
$$
\ww_{i+k/2}=A\ww_i
$$
and solve the reduced $k\times k$ eigenproblem arising by projecting onto this
$k$-dimensional subspace.
Then pick the $k/2$ most extreme eigenpairs, $(\lambda_i,\ww_i)$, $i=1,\dots,k/2$.
With these $k$ vectors, we form the reduced eigenproblem and repeat as desired.
We could of course also consider $k$ to be a multiple of $\mu$ for $\mu>2$.
We could then retain only $k/\mu$ vectors and perform $\mu-1$ Krylov steps
for each retained vector.

\subsection{Preconditioning and generalized eigenproblems}

The approach taken to $k$-step methods described via LOBPCG methods can clearly
be used to implement them.
Thus preconditioning and generalized eigenproblems can also be used.
We postpone to a later publication the study of extrapolation in these contexts.

\section{Conclusions and perspectives}

We have considered using small eigenproblems as a technique to enhance the performance
of matrix-free methods, such as the power method.
We introduce the concept of the $k$-step method, which we identify as an Arnoldi
method and also the LOBPCG method.
We have examined the performance as a function of $k$ for various test matrices.
But most importantly, we have shown that extrapolation can be used as a simple
post-processing procedure to enhance the $k$-step method at essentially zero cost.
Extrapolation can be added easily to packages such as SLEPc \cite{hernandez2005slepc}.

\section*{Acknowledgments}
SP is supported in part by the National Science Foundation NSF-DMS 1852876.


\begin{thebibliography}{10}
\bibitem{benner2011locally}
{\sc P.~Benner and T.~Mach}, {\em Locally optimal block preconditioned
  conjugate gradient method for hierarchical matrices}, PAMM, 11 (2011),
  pp.~741--742.

\bibitem{bjorck1994numerics}
{\sc {\AA}.~Bj{\"o}rck}, {\em Numerics of {Gram--Schmidt} orthogonalization},
  Linear Algebra and Its Applications, 197 (1994), pp.~297--316.

\bibitem{lrsBIBhi}
{\sc E.~Canc{\` e}s and L.~R. Scott}, {\em {van der Waals interactions between
  two hydrogen atoms: The Slater-Kirkwood} method revisited}, SIAM Journal on
  Mathematical Analysis, 50 (2018), pp.~381--410.

\bibitem{DH11}
{\sc T.~A. Davis and Y.~Hu}, {\em The University of Florida sparse matrix
  collection}, ACM Transactions on Mathematical Software (TOMS), 38 (2011),
  pp.~1--25.

\bibitem{DSHMRX19}
{\sc C.~De~Sa, B.~He, I.~Mitliagkas, C.~R{\'e}, and P.~Xu}, {\em Accelerated
  stochastic power iteration}, Proceedings of Machine Learning Research, 84
  (2019), pp.~58--67.

\bibitem{lrsBIBiu}
{\sc P.~Farrell, L.~Mitchell, L.~R. Scott, and F.~Wechsung}, {\em Robust
  discretization and multigrid solution for incompressible and nearly
  incompressible continua}, TBD, ? (2020), p.~??

\bibitem{hernandez2005slepc}
{\sc V.~Hernandez, J.~E. Roman, and V.~Vidal}, {\em {SLEPc: A} scalable and
  flexible toolkit for the solution of eigenvalue problems}, ACM Transactions
  on Mathematical Software (TOMS), 31 (2005), pp.~351--362.

\bibitem{knyazev2017recent}
{\sc A.~Knyazev}, {\em Recent implementations, applications, and extensions of
  the {Locally Optimal Block Preconditioned Conjugate Gradient method
  (LOBPCG)}}, arXiv preprint arXiv:1708.08354,  (2017).

\bibitem{knyazev2001toward}
{\sc A.~V. Knyazev}, {\em Toward the optimal preconditioned eigensolver:
  Locally optimal block preconditioned conjugate gradient method}, SIAM Journal
  on Scientific Computing, 23 (2001), pp.~517--541.

\bibitem{ref:imprestartArnoldi}
{\sc R.~B. Lehoucq and D.~C. Sorensen}, {\em Deflation techniques for an
  implicitly restarted arnoldi iteration}, SIAM Journal on Matrix Analysis and
  Applications, 17 (1996), pp.~789--821.

\bibitem{ref:Pollackextraeig}
{\sc N.~Nigam and S.~Pollock}, {\em A simple extrapolation method for clustered
  eigenvalues}, 2020, \url{https://arxiv.org/abs/2006.10164}.

\bibitem{paige2006modified}
{\sc C.~C. Paige, M.~Rozlozn{\'\i}k, and Z.~Strakos}, {\em Modified
  {Gram-Schmidt (MGS)}, least squares, and backward stability of {MGS-GMRES}},
  SIAM Journal on Matrix Analysis and Applications, 28 (2006), pp.~264--284.

\bibitem{lrsBIBjf}
{\sc S.~Pollock and L.~R. Scott}, {\em Using small eigenproblems to accelerate
  power method iterations}, Research Report {UC/CS} TR-2021, Dept. Comp.
  Sci., Univ. Chicago, 2021.

\bibitem{saad1980rates}
{\sc Y.~Saad}, {\em On the rates of convergence of the {Lanczos and the
  block-Lanczos} methods}, SIAM Journal on Numerical Analysis, 17 (1980),
  pp.~687--706.

\bibitem{saad1980arnoldi}
{\sc Y.~Saad}, {\em Variations on {Arnoldi's} method for computing
  eigenelements of large unsymmetric matrices}, Linear algebra and its
  applications, 34 (1980), pp.~269--295.

\bibitem{lrsBIBis}
{\sc L.~R. Scott}, {\em Kinetic energy flow instability with application to
  {Couette} flow}, Research Report {UC/CS} TR-2020-07, Dept. Comp. Sci., Univ.
  Chicago, 2020.

\bibitem{ref:sorensensurvey}
{\sc D.~C. Sorensen}, {\em Numerical methods for large eigenvalue problems},
  Acta Numerica, 11 (2002), p.~519.

\end{thebibliography}
\end{document}